\documentclass[a4paper,12pt]{article}

\usepackage[english]{babel} 
\usepackage[ansinew]{inputenc}
\usepackage[LGR,T1]{fontenc}

\usepackage{verbatim}
\usepackage{color}
\usepackage{calc}
\usepackage{hyperref}
\usepackage[dvips]{graphicx}
\hypersetup{colorlinks=true, linkcolor=blue, filecolor=blue, urlcolor=blue}

\usepackage{amssymb}
\usepackage[mathscr]{eucal}
\usepackage{amscd,amsmath,amsfonts,amssymb,textcomp}
\everymath{\displaystyle}

\usepackage[all]{xy}

\newcommand{\mb}[1]{\mathbb{#1}}

\newcommand{\mf}[1]{\mathfrak{#1}}

\newcommand{\MBC}{\mathbb{C}}
\newcommand{\MBF}{\mathbb{F}}

\newcommand{\MBQ}{\mathbb{Q}}

\newcommand{\MBZ}{\mathbb{Z}}

\newcommand{\MCE}{\mathcal{E}}
\newcommand{\MCF}{\mathcal{F}}

\newcommand{\MCJ}{\mathcal{J}}

\newcommand{\MCO}{\mathcal{O}}
\newcommand{\MCP}{\mathcal{P}}

\newcommand{\MCS}{\mathcal{S}}

\newcommand{\MFf}{\mathfrak{f}}

\newcommand{\MFm}{\mathfrak{m}}

\newcommand{\MFp}{\mathfrak{p}}

\newcommand{\MFq}{\mathfrak{q}}

\newcommand{\MFr}{\mathfrak{r}}

\newcommand{\MFs}{\mathfrak{s}}

\newcommand{\MSO}{\mathscr{O}}

\newcommand{\GGa}{\alpha}

\newcommand{\GVe}{\varepsilon}

\newcommand{\GGg}{\gamma}

\newcommand{\GGs}{\sigma}

\newcommand{\GGW}{\Omega}

\newcommand{\vid}{\varnothing}

\newcommand{\ho}[2]{{\mathop{#1}\limits^{#2}}}

\newcommand{\la}{\left\langle}
\newcommand{\ra}{\right\rangle}

\newcommand{\Cl}{\mathrm{Cl}}

\newcommand{\Cok}{\mathrm{Cok}}

\newcommand{\Fit}{\mathrm{Fit}}

\newcommand{\Gal}{\mathrm{Gal}}

\newcommand{\IM}{\mathrm{Im}}

\newcommand{\Ker}{\mathrm{Ker}}

\newcommand{\NN}{\mathrm{N}}

\newcommand{\PR}{\mathrm{pr}}

\usepackage{makeidx}
\title{The $\chi$-part of the analytic class number formula, for global function fields.}
\author{St\'ephane VIGUI\'E
\footnote{S.Vigui\'e, Laboratoire de math\'ematiques de Besançon, UMR CNRS 6623, Universit\'e de Franche-Comt\'e, 16 route de Gray, 25030 Besançon cedex, France.
e-mail: \texttt{stephane.viguie@univ-fcomte.fr}}
}
\makeindex

\newtheorem{dfe}{dfe}[section]
\newtheorem{lem}{lem}[section]
\newtheorem{pro}{pro}[section]
\newtheorem{rem}{rem}[section]
\newtheorem{teh}{teh}[section]

\newtheorem{df}[dfe]{Definition}
\newtheorem{lm}[lem]{Lemma}
\newtheorem{pr}[pro]{Proposition}

\newtheorem{rmq}[rem]{Remark}
\newtheorem{theo}[teh]{Theorem}

\numberwithin{equation}{section}

\begin{document}
\maketitle

\begin{abstract}
Let $F/k$ be a finite abelian extension of global function fields, totally split at a distinguished place $\infty$ of $k$.
We show that a complex Gras conjecture holds for Stark units, and we derive a refined analytic class number formula.
\end{abstract}
\bigskip

\noindent{\small\textbf{Mathematics Subject Classification (2010):} 11R58, 11R29, (11M06).}
\\

\noindent{\small\textbf{Key words:} Stark units, Gras conjecture, Analytic class number formula.}

\section{Introduction.}

Let $k$ be a global function field with constant field $\MBF_q$, and let $\infty$ be a distinguished place of $k$.
We write $k_\infty$ for the completion of $k$ at $\infty$.
For any finite abelian extension $K/k$, let $\MCO_K$ be the ring of functions of $k$ which are regular outside the places of $K$ sitting above $\infty$.
We denote by $\Cl\left(\MCO_K\right)$ the ideal-class group of $\MCO_K$.
If $K\subseteq k_\infty$ we define in section \ref{elliptic units} a group $\MCE_K$ of Stark units, which have finite index in $\MCO_K^\times$, the group of units of $\MCO_K$.
This group $\MCE_K$ has already been studied in \cite{oukhaba92}, \cite{oukhaba95}.

We fix a finite abelian extension $F\subseteq k_\infty$ of $k$, with Galois group $G$, and degree $g$.
In \cite{viguie10}, for every nontrivial irreducible rational character $\psi$ of $G$, we established that
\begin{equation}
\#\left(\MBZ_{\la g\ra}\otimes_\MBZ\left(\MCO_F^\times/\MCE_F\right)\right)_\psi = \#\left(\MBZ_{\la g\ra}\otimes_\MBZ \Cl(\MCO_F)\right)_\psi,
\label{psipart}
\end{equation}
where $\MBZ_{\la g\ra}:=\MBZ\left[g^{-1}\right]$, and where the index $\psi$ means we take the $\psi$-parts.
In \cite{oukhaba-viguie10} we used Euler systems to prove that $\MCE_F$ satisfies the Gras conjecture, 
\begin{equation*}
\#\left(\MBZ_p\otimes_\MBZ\left(\MCO_F^\times/\MCE_F\right)\right)_\psi = \#\left(\MBZ_p\otimes_\MBZ \Cl(\MCO_F)\right)_\psi,
\end{equation*}
for all prime number $p\nmid qg$ and all irreducible nontrivial $\MBQ_p$-character $\psi$
(but we were not able to prove the conjecture in the special case where the following conditions are simultaneously satisfied: 
$p\vert\#\Cl\left(\MCO_k\right)$, $\psi$ is a conjugate of the Teichmuller character, $\mu_p\not\subset k$ and $\mu_p\subset F$, 
where $\mu_p$ is the group of $p$-th roots of unity in the separable closure of $k$).

Let $\mu_g$ be the group of $g$-th roots of unity in the field of complex numbers.
Let $\MSO$ be the integral closure of $\MBZ_{\la g\ra}$ in $\mb{Q}(\mu_g)$.
For any module $M$ over a commutative ring $A$, we denote by $\Fit_A(M)$ the Fitting ideal of $M$.
Using the properties of Rubin-Stark units stated in \cite{popescu99a} and \cite{popescu99b}, we prove here a complex version of the Gras conjecture for $\MCE_F$ (Theorem \ref{theosupergras}). 
More precisely, we prove that for every nontrivial complex character $\chi$ of $G$, we have
\begin{equation}
\Fit_\MSO\left(\MSO\otimes_\MBZ\left(\MCO_F^\times/\MCE_F\right)\right)_\chi = \Fit_\MSO\left(\MSO\otimes_\MBZ \Cl(\MCO_F)\right)_\chi,
\label{chipart}
\end{equation}
which is a refinement of (\ref{psipart}).
From (\ref{chipart}), one can easily deduce that the classical Gras conjecture holds for all prime number $p\nmid g$.
Thus Theorem \ref{theosupergras} expands the result obtained in \cite{oukhaba-viguie10}, in particular it shows that the conjecture holds for $p$ equal to the caracteristic of $k$, and also for the conjugates of the Teichmuller character in the above special case.
Let us also mention the analogy between this theorem and a recent result of P.\,Buckingham (see \cite[Theorem 7.1]{buckingham11}), who is concerned with Rubin-Stark elements in cyclic extensions of totally real number fields.

Also we can combine this complex Gras conjecture (\ref{chipart}) with the computations made in \cite{viguie10}.
Thus we derive a \guillemotleft $\chi$-part\guillemotright\, version of the analytic class number formula (Theorem \ref{theosuperanclfor}), 
connecting the Fitting ideal of the $\chi$-part of the $\MSO$-module $\MSO\otimes_\MBZ \Cl\left(\MCO_F\right)$ to the value at $0$ of an $L$-function attached to $\overline{\chi}$, for any non trivial irreducible complex character $\chi$ of $G$.

\section{Stark units in function fields.}\label{elliptic units}

Let $d$ be the degree of $\infty$ over $\MBF_q$.
If $\MFm$ is a nonzero ideal of $\MCO_k$, then we denote by $H_\MFm\subseteq k_\infty$ the maximal abelian extension of $k$ contained in $k_\infty$, such that the conductor of $H_\MFm/k$ divides $\MFm$.
The function field version of the abelian conjectures of Stark, proved by P.\,Deligne in \cite{tate84} by using \'etale cohomology or by D.\,Hayes in \cite{hayes85} by using Drinfel'd modules, 
claims that, for any proper nonzero ideal $\MFm$ of $\MCO_k$, there exists an element $\GVe_\MFm\in H_\MFm$, unique up to roots of unity, such that
\\ 

\noindent $\left.i\right)$ The extension $H_\MFm\left(\GVe_\MFm^{1/w_\infty}\right)/k$ is abelian, where $w_\infty:=q^d-1$.
Moreover, it is unramified outside $S_\MFm$, where $S_\MFm$ is the set containing $\infty$ and the places of $k$ which divide $\MFm$.
\\

\noindent $\left.ii\right)$ If $\MFm$ is divisible by two prime ideals then $\GVe_\MFm$ is a unit of $\MCO_{H_\MFm}$.
If $\MFm=\MFq^e$, where $\MFq$ is a prime ideal of $\MCO_k$ and $e$ is a positive integer, then
\[\GVe_\MFm\MCO_{H_\MFm}=(\MFq)_\MFm^{\frac{w_\infty}{w_k}},\]
where $w_k:=q-1$ and $(\MFq)_\MFm$ is the product of the prime ideals of $\MCO_{H_\MFm}$ which divide $\MFq$.
\\

\noindent $\left.iii\right)$ We have
\begin{equation} 
\mathrm{L}_\MFm(0,\chi)=\frac{1}{w_\infty}\sum_{\GGs\in \Gal(H_\MFm/k)}\chi(\GGs)v_\infty\left(\GVe_\MFm^\GGs\right) 
\label{equaLMFmGVEMFm} 
\end{equation}
for all complex irreducible characters of $\Gal\left(H_\MFm/k\right)$, where $v_\infty$ is the normalized valuation of $k_\infty$.
\\

Let us recall that $s\mapsto \mathrm{L}_\MFm(s,\chi)$ is the $L$-function associated to $\chi$, defined for the complex numbers $s$ such that $\mathrm{Re}(s)>1$ by the Euler product
\[\mathrm{L}_\MFm(s,\chi)=\prod_{v\nmid\MFm}\left(1-\chi(\GGs_v)\NN(v)^{-s}\right)^{-1},\]
where $v$ describes the set of places of $k$ not dividing $\MFm$.
For such a place, $\GGs_v$ and $\NN(v)$ are the Frobenius automorphism of $H_\MFm/k$ and the order of the residue field at $v$ respectively.
Let us remark that $\GGs_\infty=1$ and $\NN(\infty)=q^d$.

For any finite abelian extension $L$ of $k$ we denote by $\MCJ_L\subseteq\mb{Z}\left[\Gal(L/k)\right]$ the annihilator of $\mu(L)$, the group of roots of unity in $L$.
The description of $\MCJ_L$ given in \cite[Lemma 2.5]{hayes85} and the property $\left.i\right)$ of $\GVe_\MFm$ implies that for any $\eta\in\MCJ_{H_\MFm}$ there exists $\GVe_\MFm(\eta)\in H_\MFm$ such that
\begin{equation}
\label{bYvhrD}
\GVe_\MFm(\eta)^{w_\infty}=\GVe_\MFm^\eta.
\end{equation}

\begin{df}
Let ${\MCP_F}$ be the subgroup of $F^\times$ generated by ${\mu}(F)$ and by all the norms 
\[\GVe_{F,\MFm}(\eta) := \NN_{H_\MFm/H_\MFm\cap F}\left(\GVe_\MFm(\eta)\right),\] 
where $\MFm$ is any nonzero proper ideal of $\MCO_{k}$, and $\eta\in\MCJ_{H_\MFm}$.
Then we set
\[{\MCE_F}:={\MCP_F}\cap \MCO_F^\times.\]
\end{df}

\section{Preliminary lemmas.}

For any finite group $H$, let us denote by $\widehat{H}$ the group of complex irreducible characters of $H$.
Then for every $\chi\in\widehat{H}$ we set $e_\chi:=\frac{1}{\#H}\sum_{\GGs\in G}\chi(\GGs)\GGs^{-1}$.
In case $H=G$, then $e_\chi$ belongs to $\MSO[G]$.
Moreover, if $\zeta\in\mu_g$ is such that $\zeta\neq1$, then $\left(1-\zeta\right)\in\MSO^\times$, thanks to the formula $g=\prod_{\substack{\zeta\in\mu_g\\\zeta\neq1}}(1-\zeta)$.
For any $\MSO[G]$-module $M$ let $M_\chi:=e_\chi M$ be the $\chi$-part of $M$.

For $K/k$ a finite abelian extension, and $S\neq\vid$ a finite set of places of $k$, let $\MCO_{K,S}$ be the Dedekind ring of the $S$-integers of $K$, i.e the functions $f\in K$ which only poles are at the places sitting over $S$.
Let $\Cl\left(\MCO_{K,S}\right)$ be the ideal class group of $\MCO_{K,S}$.

\begin{lm}\label{lmSSMFmcepareil}
Let $\MFm$ be a nonzero ideal of $\MCO_k$, and let $\chi\in\widehat{G}$, $\chi\neq1$.
Assume that for every prime ideal $\MFp$ of $\MCO_k$ which divides $\MFm$, $\chi$ is not trivial on the decomposition group $D_\MFp$ of $\MFp$ in $F/k$.
Let $S_\MFm$ be the set of places of $k$ which contains $\infty$ and all the prime divisors of $\MFm$.
Then in the category of $\MSO[G]$-modules, we have
\[\left(\MSO\otimes_\MBZ\MCO_{F,S_\MFm}^\times\right)_\chi = \left(\MSO\otimes_\MBZ\MCO_F^\times\right)_\chi \quad\text{and}\quad \left(\MSO\otimes_\MBZ \Cl\left(\MCO_{F,S_\MFm}\right)\right)_\chi \simeq \left(\MSO\otimes_\MBZ \Cl\left(\MCO_F\right)\right)_\chi.\]
\end{lm}

\noindent\textsl{Proof.}
Let $\MCS$ be the $\MBZ[G]$-module generated by the prime ideals of $\MCO_F$ dividing $\MFm\MCO_F$.
Let $\bar{\MCS}$ be the image of $\MCS$ in $\Cl\left(\MCO_F\right)$.
We have the following exact sequences:
\[\xymatrix@!0 @R=6mm @C=1cm {
0 \ar[rr] && \bar{\MCS} \ar[rrr] &&& \Cl\left(\MCO_F\right) \ar[rrr] &&& \Cl\left(\MCO_{F,S_\MFm}\right) \ar[rr] && 0, \\
0 \ar[rr] && \MCO_F^\times \ar[rrr] &&& \MCO_{F,S_\MFm}^\times \ar[rrr] &&& \MCS, && \\
 && &&& x \ar@{|->}[rrr] &&& \prod_{\MFp\in\MCS}\MFp^{v_\MFp(x)}, && 
}\]
where for all $\MFp\in\MCS$, $v_\MFp$ is the normalized valuation at $\MFp$.
Since $\MSO$ is $\MBZ$-flat, all we have to show is $\left(\MSO\otimes_\MBZ\MCS\right)_\chi=0$.
Let $\MFp\in\MCS$.
There is $\GGg\in G$ such that $\MFp^\GGg=\MFp$ and $\chi(\GGg)\neq1$.
Then $(\chi(\GGg)-1)e_\chi(1\otimes\MFp)= e_\chi(1\otimes\MFp^{\GGg-1}) =0$, with $\left(\chi(\GGg)-1\right)\in\MSO^\times$, hence $e_\chi\left(1\otimes\MFp\right) = 0$.
\hfill$\square$\\

For $K\subseteq k_\infty$ a finite abelian extension of $k$, let $\ell_K:K^\times\rightarrow\MBZ\left[\Gal(K/k)\right]$ be the $\Gal(K/k)$-equivariant map defined by 
\[\ell_K(x):=\sum_{\GGs\in \Gal(K/k)}v_\infty\left(x^\GGs\right)\GGs^{-1}.\]

\begin{lm}\label{lmstarkamochee}
Let $\MFm$ be a nonzero ideal of $\MCO_k$, and $\GGa\in\MCJ_{H_\MFm}$.
We set $C_\MFm:=\Cl \left(\MCO_{F\cap H_\MFm,S_\MFm}\right)$.
Then 
\[\ell_F\left(\GVe_{F,\MFm}(\GGa)\right) \in \Fit_{\MBZ_{\la g\ra}[G]}(C_\MFm) \ell_F\left(\MCO_{F\cap H_\MFm,S_\MFm}^\times\right).\]
\end{lm}

\noindent\textsl{Proof.}
By the description of $\MCJ_{H_\MFm}$ given in \cite{hayes85}, Lemma 2.5, we know that there are a finite set $T$ of nonzero prime ideals of $\MCO_k$, and a family $\left(\GGa_\MFp\right)_{\MFp\in T}\in\MBZ\left[\Gal\left(H_\MFm/k\right)\right]^T$, 
such that $S_\MFm\cap T=\vid$ and $\GGa = \sum_{\MFp\in T}\GGa_\MFp\left(1-\NN(\MFp)\GGs_\MFp^{-1}\right)$,
where $\GGs_\MFp$ is the Frobenius of $\MFp$ in $H_\MFm/k$.
It suffices to show $\ell_F\left(\GVe_{F,\MFm}\left(1-\NN(\MFp)\GGs_\MFp^{-1}\right)\right) \in \Fit_{\MBZ_{\la g\ra}[G]}(C_\MFm) \ell_F\left(\MCO_{F\cap H_\MFm,S_\MFm}^\times\right)$, for a fixed nonzero prime ideal $\MFp$ of $\MCO_k$, with $\MFp\notin S_\MFm$.

For any abelian extension $K$ of $k$, we denote by $U_{K,\MFm,\MFp}$ the group of units of $\MCO_{K,S_\MFm}$ which are congruent to $1$, modulo all primes above $\MFp$.
From (\ref{bYvhrD}) we deduce that $\GVe_\MFm\left(1-\NN(\MFp)\GGs_\MFp^{-1}\right)$ can be chosen such that $\GVe_\MFm\left(1-\NN(\MFp)\GGs_\MFp^{-1}\right) \in U_{H_{\MFm},\MFm,\MFp}$.
For $\chi\in\widehat{\Gal\left(H_\MFm/k\right)}$, we define the meromorphic function $s\mapsto \mathrm{L}_{\MFm,\MFp}(s,\overline{\chi})$ by \[\mathrm{L}_{\MFm,\MFp}(s,\overline{\chi}) := \left(1-\NN(\infty)^{-s}\right) \mathrm{L}_\MFm(s,\overline{\chi}) \left(1-\NN(\MFp)^{1-s}\overline{\chi}\left(\sigma_\MFp\right)\right).\]

Derivating, and using the property $\left.iii\right)$ of Stark units, we obtain
\begin{eqnarray}
\mathrm{L}'_{\MFm,\MFp}(0,\overline{\chi}) & = & d.\mathrm{ln}(q)\mathrm{L}_\MFm(0,\overline{\chi}) \left(1-\NN(\MFp)\chi\left(\sigma_\MFp^{-1}\right)\right)\nonumber\\
& = & d.\mathrm{ln}(q) \chi\left(\ell_{H_{\MFm}}\left(\GVe_{\MFm}\left(1-\NN(\MFp)\GGs_\MFp^{-1}\right)\right)\right),
\label{equarubstarksol}
\end{eqnarray}
where $\chi$ is extended to $\MBC\left[\Gal(H_\MFm/k)\right]$ by linearity.
For $s\in\MBC$, we set 
\[\Theta_{\MFm,\MFp}(s) = \sum_{\chi\in\widehat{\Gal\left(H_\MFm/k\right)}}\mathrm{L}_{\MFm,\MFp}(s,\overline{\chi})e_\chi \quad\text{in}\quad \MBC\left[\Gal\left(H_\MFm/k\right)\right],\] 
wherever it is defined.
From (\ref{equarubstarksol}), we have
\begin{eqnarray*}
\Theta_{\MFm,\MFp}'(0) & = & d.\mathrm{ln}(q)\ell_{H_{\MFm}}\left(\GVe_{\MFm}\left(1-\NN(\MFp)\GGs_\MFp^{-1}\right)\right)\\
& = & -\sum_{\GGg\in \Gal\left(H_{\MFm}/{k}\right)} \mathrm{ln} \left(\left| \GVe_{\MFm}\left(1-\NN(\MFp)\GGs_\MFp^{-1}\right)^{\GGg^{-1}} \right|_{\infty}\right)\GGg.
\end{eqnarray*}
We have just verified that $\GVe_{\MFm}\left(1-\NN(\MFp)\GGs_\MFp^{-1}\right)$ satisfies \cite[Theorem 0]{popescu99b}, for $\left(H_\MFm,S_\MFm,\{\MFp\},1\right)$. 
Then $\GVe_{F,\MFm}\left(1-\NN(\MFp)\GGs_\MFp^{-1}\right)$ satisfies \cite[Theorem 0]{popescu99b}, for $\left(H_\MFm\cap F,S_\MFm,\{\MFp\},1\right)$, and we have
\[\ell_{F}\left(\GVe_{F,\MFm}\left(1-\NN(\MFp)\GGs_\MFp^{-1}\right)\right) \in \Fit_{\MBZ_{\la g\ra}[G]} \left(\MBZ_{\la g\ra}\otimes_\MBZ \Cl \left(\MCO_{F\cap H_\MFm,S_\MFm},\MFp\right)\right) \ell_F\left(U_{F\cap H_\MFm,\MFm,\MFp}^\times\right),\]
where $\Cl \left(\MCO_{F\cap H_\MFm,S_\MFm},\MFp\right)$ is the quotient of fractional ideals of $\MCO_{F\cap H_\MFm}$ by those principal fractional ideals, which are generated by an element congruent to $1$ modulo $\MFp$.
But $\MBZ_{\la g\ra}\otimes_\MBZ C_\MFm$ is a quotient of $\MBZ_{\la g\ra}\otimes_\MBZ \Cl \left(\MCO_{F\cap H_\MFm,S_\MFm},\MFp\right)$, and $\MBZ_{\la g\ra}[G]$ is a finite product of Dedekind rings, so 
\[\Fit_{\MBZ_{\la g\ra}[G]} \left(\MBZ_{\la g\ra}\otimes_\MBZ \Cl \left(\MCO_{F\cap H_\MFm,S_\MFm},\MFp\right)\right) \subseteq \Fit_{\MBZ_{\la g\ra}[G]} \left(\MBZ_{\la g\ra}\otimes_\MBZ C_\MFm\right),\]
and we are done.
\hfill$\square$\\

\begin{df}
Let $\GGW$ be the $\MBZ[G]$-submodule of $F^\times$ generated by $\mu(F)$ and by the elements $\GVe_{F,\MFm}:=\NN_{H_\MFm/F\cap H_\MFm}(\GVe_\MFm)$, where $\MFm$ is any nonzero proper ideal of $\MCO_k$.
\end{df}

For any $\chi\in\widehat{G}$, let $F_\chi$ be the subfield of $F$ fixed by $\Ker(\chi)$, and $\MFf_\chi$ be the conductor of $F_\chi$.
$\chi_{\PR}$ denotes the character of $\Gal\left(H_{\MFf_\chi}/k\right)$ defined by $\chi$.
Assume that $\chi$ is nontrivial.
Then by \cite[Proposition 3.1]{viguie10}, we have
\begin{equation}
\label{proptrun}
\left(\MSO\ell_F\left(\GGW\right)\right)_\chi = \MSO w_\infty\mathrm{L}_{\MFf_\chi}\left(0,\bar{\chi}_\PR\right)e_\chi.
\end{equation}
Moreover one can easily relate $\ell_F\left(\GGW\right)$ to $\ell_F\left(\MCE_F\right)$ (see \cite[(3.13)]{viguie10}).
If $\chi\neq1$, we have
\begin{equation}
\label{GGWMCEF}
\left(\MSO\ell_F\left(\MCE_F\right)\right)_\chi = \left(\MSO w_\infty^{-1}\MCJ_F\ell_F\left(\GGW\right)\right)_\chi.
\end{equation}

\begin{lm}\label{lmgrasun}
Let $\chi\in\widehat{G}$, $\chi\neq1$.
There is a nonzero ideal $\MFm$ of $\MCO_k$, satisfying the following properties:

$\left.i\right)$ $\MFm$ is divisible by at least two distinct prime ideals, 

$\left.ii\right)$ $F_\chi \subseteq H_\MFm$, 

$\left.iii\right)$ for each prime ideal $\MFp$ which divides $\MFm$, $\chi$ is not trivial on the decomposition group $D_\MFp$ of $\MFp$ in $F/k$.

$\left.iv\right)$ As an $\MSO$-module, $\left(\MSO\otimes_\MBZ\MCE_F\right)_\chi$ is generated by $\left(\MSO\otimes_\MBZ\mu(F)\right)_\chi$ and by the elements $e_\chi\left(1\otimes\GVe_{F,\MFm}(\GGa)\right)$, where $\GGa \in \MCJ_{H_\MFm}$.
\end{lm}

\noindent\textsl{Proof.}
Since $\chi\neq1$, we can find two distinct nonzero prime ideals $\MFp$ and $\MFq$ of $\MCO_k$, unramified in $F/k$, such that 
\[\bar{\chi}_{\PR}\left(\GGs_\MFp\right)\neq 1\quad \text{and} \quad \bar{\chi}_{\PR}\left(\GGs_\MFq\right)\neq 1.\]
We set $\MFm=\MFf_\chi\MFp\MFq$.
Obviously, conditions $\left.i\right)$, $\left.ii\right)$ and $\left.iii\right)$ are satisfied.
As in the proof of \cite[Proposition 3.1]{viguie10}, from the property $\left.iii\right)$ of Stark units we obtain
\[\ell_F\left(\GVe_{F,\MFm}\right)e_\chi = [F:F\cap H_{\MFm}]w_\infty \mathrm{L}_{\MFf_\chi}(0,\overline{\chi}_{\PR})  \left(1-\bar{\chi}_{\PR}\left(\GGs_\MFp\right)\right) \left(1-\bar{\chi}_{\PR}\left(\GGs_\MFq\right)\right) e_\chi.\]
Since $[F:F\cap H_{\mf{n}\MFp}]$, $1-\bar{\chi}_{\PR}\left(\GGs_\MFp\right)$ and $1-\bar{\chi}_{\PR}\left(\GGs_\MFp\right)$ belong to $\MSO^\times$, and by (\ref{proptrun}), we have
\begin{equation}
\MSO\ell_F\left(\GVe_{F,\MFm}\right) e_\chi = \MSO w_\infty \mathrm{L}_{\MFf_\chi}(0,\overline{\chi}_{\PR}) e_\chi = \MSO\ell_F(\GGW)e_\chi.
\label{ggogoltha}
\end{equation}
From (\ref{ggogoltha}) and (\ref{GGWMCEF}), we deduce
\[\MSO\MCJ_Fw_\infty^{-1}\ell_F\left(\GVe_{F,\MFm}\right) e_\chi = \MSO\MCJ_Fw_\infty^{-1}\ell_F(\GGW)e_\chi = \MSO\ell_F(\MCE_F)e_\chi,\]
and condition $\left.iv\right)$ follows.
\hfill$\square$\\

To go further we need some preliminary remarks.

\begin{rmq}\label{fittOGO}
For any $\MSO[G]$-module $M$, we have $\Fit_{\MSO[G]}(M)=\sum_{\chi\in\widehat{G}}\Fit_\MSO(M_\chi)e_\chi$.
\end{rmq}

\begin{rmq}\label{rmqPsiHsubgrupG}
Let $H$ be a sub-group of $G$.
Let $M$ and $N$ be two $G$-modules, and $\psi:M\rightarrow N$ be a $G$-equivariant map.
If $\Cok(\Psi):=N/\IM(\Psi)$ is annihilated by $\#(H)$ then we derive from $\Psi$ a surjective map
\[\Psi_\MSO:\MSO\otimes_\MBZ M\twoheadrightarrow\MSO\otimes_\MBZ N.\]
Let us assume, in addition, that $\Ker(\Psi)$ is annihilated by $\Sigma\GGs$, $\GGs\in H$.
Then, for every $\chi\in\widehat{G}$ trivial on $H$, the restriction of $\Psi_\MSO$ gives an isomorphism
\[\left(\MSO\otimes_\MBZ M\right)_\chi \simeq \left(\MSO\otimes_\MBZ N\right)_\chi.\] 
\end{rmq}

As a particular case, for any subextension $K/k$ of $F/k$ and $H=\Gal(F/K)$, and any nonzero ideal $\MFm$ of $\MCO_k$, the norm maps give isomorphisms 
\[\left(\MSO\otimes_\MBZ \Cl(\MCO_{F,S_\MFm})\right)_\chi \ho{\longrightarrow}{\sim} \left(\MSO\otimes_\MBZ \Cl(\MCO_{K,S_\MFm})\right)_\chi \quad\text{and}\quad \left(\MSO\otimes_\MBZ \MCO_{F,S_\MFm}^\times\right)_\chi \ho{\longrightarrow}{\sim} \left(\MSO\otimes_\MBZ \MCO_{K,S_\MFm}^\times\right)_\chi,\]
for any $\chi\in\widehat{G}$ which is trivial on $H$.
Since $\#(H)\in\MSO^\times$ we deduce that for such a character, the canonical inclusion also gives an equality
\[\left(\MSO\otimes_\MBZ \MCO_{K,S_\MFm}^\times\right)_\chi = \left(\MSO\otimes_\MBZ \MCO_{F,S_\MFm}^\times\right)_\chi.\]

\begin{pr}\label{lmgrasdeux}
Let $\chi\in\widehat{G}$, $\chi\neq1$.
We have 
\[\left(\MSO\ell_F\left(\MCE_F\right)\right)_\chi \subseteq \Fit_\MSO\left(\left(\MSO\otimes_\MBZ \Cl\left(\MCO_F\right)\right)_\chi\right) \left(\MSO\ell_F\left(\MCO^\times_F\right)\right)_\chi.\]
\end{pr}

\noindent\textsl{Proof.}
We choose an ideal $\MFm$ of $\MCO_k$ satisfying the four conditions of Lemma \ref{lmgrasun}.
Because of Lemma \ref{lmgrasun}, $\left.iv\right)$, it is sufficient to show that for $\GGa \in \MCJ_{H_\MFm}$, we have 
\[\ell_F(\GVe_{F,\MFm}(\GGa))e_\chi \in \Fit_\MSO\left(\left(\MSO\otimes_\MBZ \Cl\left(\MCO_F\right)\right)_\chi\right) \left(\ell_F\left(\MCO^\times_F\right)\right)_\chi.\]
By Lemma \ref{lmstarkamochee}, we know that 
\[\ell_F(\GVe_{F,\MFm}(\GGa)) \in \Fit_{\MBZ[G]}\left(C_\MFm\right) \ell_F\left(\MCO_{F\cap H_\MFm,S_\MFm}^\times\right).\] 
From Remark \ref{fittOGO}, we deduce 
\[\ell_F(\GVe_{F,\MFm}(\GGa))e_\chi \in \Fit_{\MSO}\left(\left(\MSO\otimes_\MBZ C_\MFm\right)_\chi\right) \left(\MSO\ell_F\left(\MCO_{F\cap H_\MFm,S_\MFm}^\times\right)\right)_\chi.\]
By Lemma \ref{lmgrasun}, $\left.ii\right)$, and Remark \ref{rmqPsiHsubgrupG}, the norm map defines an isomorphism
\[\left(\MSO\otimes_\MBZ \Cl\left(\MCO_{F,S_\MFm}\right)\right)_\chi \simeq \left(\MSO\otimes_\MBZ C_\MFm\right)_\chi,\]
and the canonical inclusion $F\cap H_\MFm\hookrightarrow F$ gives 
\[\left(\MSO\otimes_\MBZ\MCO_{F\cap H_\MFm,S_\MFm}^\times\right)_\chi = \left(\MSO\otimes_\MBZ\MCO_{F,S_\MFm}^\times\right)_\chi.\]
Then 
\[\ell_F(\GVe_{F,\MFm}(\GGa))e_\chi \in \Fit_{\MSO}\left(\left(\MSO\otimes_\MBZ \Cl\left(\MCO_{F,S_\MFm}\right)\right)_\chi\right) \left(\MSO\ell_F\left(\MCO_{F,S_\MFm}^\times\right)\right)_\chi.\]
From Lemma \ref{lmSSMFmcepareil}, and Lemma \ref{lmgrasun}, $\left.iii\right)$, we know that
\[\Fit_{\MSO}\left(\left(\MSO\otimes_\MBZ \Cl\left(\MCO_{F,S_\MFm}\right)\right)_\chi\right) = \Fit_{\MSO}\left(\left(\MSO\otimes_\MBZ \Cl\left(\MCO_F\right)\right)_\chi\right),\]
\[\text{and } \left(\MSO\otimes_\MBZ\MCO_{F,S_\MFm}^\times\right)_\chi =  \left(\MSO\otimes_\MBZ\MCO_F^\times\right)_\chi.\]
The proposition follows.
\hfill$\square$\\

\section{Statement and proof of the theorems.}

\begin{theo}\label{theosupergras}
Let $\chi\in\widehat{G}$, $\chi\neq1$.
We have 
\[\Fit_\MSO\left(\left(\MSO\otimes_\MBZ \Cl\left(\MCO_F\right)\right)_\chi\right) = \Fit_\MSO\left( \left(\MSO\otimes_\MBZ\left(\MCO_F^\times/\MCE_F\right)\right)_\chi \right).\]
\end{theo}

\noindent\textsl{Proof.}
We have
\begin{eqnarray}
\prod_{\substack{\xi\in\widehat{G}\\\xi\neq1}} \Fit_\MSO\left(\left(\MSO\otimes_\MBZ \Cl\left(\MCO_F\right)\right)_\xi\right) & = & \Fit_\MSO\left((1-e_1).\MSO\otimes_\MBZ \Cl\left(\MCO_F\right)\right) \nonumber\\
 & = & \MSO\#\left((1-e_1).\MBZ_{\la g\ra}\otimes_\MBZ \Cl\left(\MCO_F\right)\right).
\label{prodfittcl}
\end{eqnarray}
In the same way, we have
\begin{eqnarray}
\prod_{\substack{\xi\in\widehat{G}\\\xi\neq1}} \Fit_\MSO\left(\left(\MSO\otimes_\MBZ \left(\MCO_F^\times/\MCE_F\right)\right)_\xi\right) & = & \MSO \Fit_{\MBZ_{\la g\ra}}\left(\MBZ_{\la g\ra}\otimes_\MBZ \left(\MCO_F^\times/\MCE_F\right)\right) \nonumber\\
 & = & \MSO \left[\MBZ_{\la g\ra}\otimes_\MBZ\MCO_F^\times : \MBZ_{\la g\ra}\otimes_\MBZ\MCE_F\right].
\label{prodfittmodinc}
\end{eqnarray}
By (\ref{psipart}), we know that
\begin{equation}
\#\left((1-e_1).\MBZ_{\la g\ra}\otimes_\MBZ \Cl\left(\MCO_F\right)\right) = \left[\MBZ_{\la g\ra}\otimes_\MBZ\MCO_F^\times : \MBZ_{\la g\ra}\otimes_\MBZ\MCE_F\right].
\label{clOFEFcard}
\end{equation}
Putting (\ref{prodfittcl}), (\ref{clOFEFcard}) and (\ref{prodfittmodinc}) together, we obtain
\[\prod_{\substack{\xi\in\widehat{G}\\\xi\neq1}} \Fit_\MSO\left(\left(\MSO\otimes_\MBZ \Cl\left(\MCO_F\right)\right)_\xi\right) = \prod_{\substack{\xi\in\widehat{G}\\\xi\neq1}} \Fit_\MSO\left(\left(\MSO\otimes_\MBZ \left(\MCO_F^\times/\MCE_F\right)\right)_\xi\right).\]
From this equality, it follows that it is sufficient to show the divisibility 
\begin{equation}
\Fit_\MSO\left(\left(\MSO\otimes_\MBZ \Cl\left(\MCO_F\right)\right)_\xi\right) \big| \Fit_\MSO\left( \left(\MSO\otimes_\MBZ\left(\MCO_F^\times/\MCE_F\right)\right)_\xi \right),
\label{divisibility}
\end{equation}
for all $\xi\in\widehat{G}$.
Since $\MCO_F^\times\cap\Ker(\ell_F)=\mu(F)$, we have $\MCO_F^\times/\MCE_F \simeq \ell_F\left(\MCO_F^\times\right) / \ell_F\left(\MCE_F\right)$ and
\begin{equation}
\label{bonobo}
\left(\MSO\otimes_\MBZ\left(\MCO_F^\times/\MCE_F\right)\right)_\xi \simeq \left(\MSO\ell_F\left(\MCO_F^\times\right)\right)_\xi / \left(\MSO\ell_F\left(\MCE_F\right)\right)_\xi.
\end{equation}
We set $\MCF_\xi:=\Fit_\MSO\left(\left(\MSO\otimes_\MBZ \Cl\left(\MCO_F\right)\right)_\xi\right)$ for convenience.
From Proposition \ref{lmgrasdeux} we derive the tautological exact sequence
\begin{equation}
\label{gifrgae}
0 \longrightarrow \frac{\MCF_\xi\left(\MSO\ell_F\left(\MCO_F^\times\right)\right)_\xi} {\left(\MSO\ell_F\left(\MCE_F\right)\right)_\xi} \longrightarrow \frac{\left(\MSO\ell_F\left(\MCO_F^\times\right)\right)_\xi} {\left(\MSO\ell_F\left(\MCE_F\right)\right)_\xi} \longrightarrow \frac{\left(\MSO\ell_F\left(\MCO_F^\times\right)\right)_\xi} {\MCF_\xi\left(\MSO\ell_F\left(\MCO_F^\times\right)\right)_\xi} \longrightarrow 0.
\end{equation}
Since $\MSO$ is a Dedekind ring, we deduce from (\ref{gifrgae}) and (\ref{bonobo}) that
\begin{eqnarray}
\label{OnUYb}
\Fit_\MSO\left(\left(\MSO\otimes_\MBZ\left(\MCO_F^\times/\MCE_F\right)\right)_\xi\right) & = &
\Fit_\MSO\left(\frac{\MCF_\xi\left(\MSO\ell_F\left(\MCO_F^\times\right)\right)_\xi} {\left(\MSO\ell_F\left(\MCE_F\right)\right)_\xi}\right) 
\Fit_\MSO\left(\frac{\left(\MSO\ell_F\left(\MCO_F^\times\right)\right)_\xi} {\MCF_\xi\left(\MSO\ell_F\left(\MCO_F^\times\right)\right)_\xi}\right) 
\nonumber\\
& = & \Fit_\MSO\left(\frac{\MCF_\xi\left(\MSO\ell_F\left(\MCO_F^\times\right)\right)_\xi} {\left(\MSO\ell_F\left(\MCE_F\right)\right)_\xi}\right) 
\MCF_\xi.
\end{eqnarray}
\hfill$\square$\\

Before we go further, we state here the definition and basic properties of index-modules, which we introduced in \cite{viguie10}.
We refer the reader to \cite{viguie10} for the proofs.

Let $K$ be a commutative field, $A\subseteq K$ be a Dedekind ring and $V$ be a $K$-vector space.
By an $A$-lattice of $V$, we mean a finitely generated $A$-submodule $R$ of $V$, such that the $K$-vector subspace of $V$ generated by $R$, denoted by $KR$, has dimension equal to the $A$-rank of $R$.

\begin{df}
Let $R\neq0$ and $S$ be $A$-lattices of $V$.
We call $A$-index-module of the couple $(R,S)$ the set
\[[R:S]_A\,\,:=\,\,\{\det(u);u\in End_{K}(V')/u(R)\subseteq S\},\]
where $V'$ is the $K$-subspace of $V$ generated by $R$ and $S$, $V'=KR+KS$.
\end{df}

In fact, $[R:S]_A$ is an $A$-submodule of $K$, and we have the following properties, for any $A$-lattices $R$, $S$, and $T$ of $V$, with $R\neq0$.
\\

$\left.i\right)$  $[R:R]_A=A$.
\\

$\left.ii\right)$  If $KS\subseteq KR$, then $[R:S]_A$ is a finitely generated $A$-submodule of $K$. 
Moreover, its $A$-rank is $1$ if $KR=KS$, and $[R:S]_A=0$ if $KS\subsetneq KR$.
\\

$\left.iii\right)$  Assume that $KR=KS$, and that $\MFr$ and $\MFs$ are two nonzero fractional ideals of $A$.
Then $[\MFr R:\MFs S]_A = \MFs^d\MFr^{-d}[R:S]_A$, where $d$ is the common $A$-rank of $R$ and $S$.
\\

$\left.iv\right)$ If $S\neq0$, and $KT\subseteq KS\subseteq KR$, then $[R:T]_A=[R:S]_A[S:T]_A$.
\\

$\left.v\right)$ If $S\subseteq R$, then $[R:S]_A=\Fit_A(R/S)$.
\\

In the sequel, we are concerned with the following situation.
$A:=\MSO$, $V:=\MBC[G]$, and the $\MSO$-lattices are $\left(\MSO\ell_F(\MCO_F^\times)\right)_\chi$,  $\left(\MSO\ell_F(\MCE_F)\right)_\chi$, $\left(\MSO\ell_F(\GGW)\right)_\chi$, ..., where $\chi$ is a nontrivial complex character of $G$.
They are all of $\MSO$-rank $1$.\\

\begin{lm}\label{lmsuperanclfor}
Let $\chi\in\widehat{G}$, $\chi\neq1$.
We have 
\[\left[\left(\MSO \ell\left(\GGW\right)\right)_\chi :  \left(\MSO\ell\left(\MCE_F\right)\right)_\chi \right]_\MSO = w_\infty^{-1} \Fit_\MSO\left(\left(\MSO\otimes_\MBZ\mu(F)\right)_\chi\right),\]
where $\left[\left(\MSO \ell\left(\GGW\right)\right)_\chi :  \left(\MSO\ell\left(\MCE_F\right)\right)_\chi \right]_\MSO$ and $w_\infty^{-1} \Fit_\MSO\left(\left(\MSO\otimes_\MBZ\mu(F)\right)_\chi\right)$ are viewed as $\MSO$-submodules of $\MBQ(\mu_g)$.
\end{lm}

\noindent\textsl{Proof.}
Let $\zeta$ be a primitive $w_F$-th root of unity in $F$, where $w_F:=\#\left(\mu(F)\right)$.
We have an obvious exact sequence,
\[\xymatrix @!0 @R=6mm @C=1cm {
0 \ar[rr] && \MSO\MCJ_F e_\chi \ar[rrr] &&& \MSO e_\chi \ar[rrr] &&& \left(\MSO\otimes_\MBZ\mu(F)\right)_\chi \ar[rr] && 0, \\
 && &&& \GGa e_\chi \ar@{|->}[rrr] &&& e_\chi\left(\GGa\otimes\zeta\right).  && \\
}\]
Using the property $\left.v\right)$ of index-modules, we deduce 
\begin{equation}
\Fit_\MSO\left(\left(\MSO\otimes_\MBZ\mu(F)\right)_\chi\right) = \left[\MSO e_\chi:\MSO\MCJ_F e_\chi\right]_\MSO.
\label{equafittmu}
\end{equation}
Also, using (\ref{GGWMCEF}), we have
\[\left[\left(\MSO \ell_F\left(\GGW\right)\right)_\chi :  \left(\MSO\ell_F\left(\MCE_F\right)\right)_\chi \right]_\MSO = \left[\left(\MSO \ell_F\left(\GGW\right)\right)_\chi :  \left(\MSO\MCJ_Fw_\infty^{-1}\ell_F\left(\GGW\right)\right)_\chi \right]_\MSO.\]
By (\ref{proptrun}), the property $\left.iii\right)$ of index-modules, and (\ref{equafittmu}) for the last equality, we deduce
\begin{eqnarray*}
\left[\left(\MSO \ell_F\left(\GGW\right)\right)_\chi :  \left(\MSO\ell_F\left(\MCE_F\right)\right)_\chi \right]_\MSO & = & \left[\MSO w_\infty \mathrm{L}_{\MFf_\chi}(0,\overline{\chi}_{\PR})e_\chi :  \MSO\MCJ_F\mathrm{L}_{\MFf_\chi}(0,\overline{\chi}_{\PR})e_\chi \right]_\MSO \\
& = & w_\infty^{-1}\left[\MSO e_\chi : \MSO\MCJ_F e_\chi \right]_\MSO \\
& = & w_\infty^{-1} \Fit_\MSO\left(\left(\MSO\otimes_\MBZ\mu(F)\right)_\chi\right).
\end{eqnarray*}
\hfill$\square$\\

The regulator $\mathrm{R}(\MCO_F)$ of $\MCO_F$ is known to be equal to $\left[\MBZ[G]_0:\ell_F\left(\MCO_F^\times\right)\right]$, where $\MBZ[G]_0$ is the augmentation ideal of $\MBZ[G]$.
Hence it is natural to take $\mathrm{R}(\MCO_F)_\chi := \left[\MSO e_\chi : \ell_F\left(\MCO_F^\times\right)e_\chi\right]_\MSO$ as a definition for the \guillemotleft $\chi$-part\guillemotright\, of the regulator, for any nontrivial character $\chi\in\widehat{G}$.

\begin{theo}\label{theosuperanclfor}
Let $\chi\in\widehat{G}$, $\chi\neq1$.
Then we have 
\[\mathrm{R}(\MCO_F)_\chi\Fit_\MSO\left(\left(\MSO\otimes_\MBZ \Cl\left(\MCO_F\right)\right)_\chi\right) =  \Fit_\MSO\left(\left(\MSO\otimes_\MBZ\mu(F)\right)_\chi\right) \mathrm{L}_{\MFf_\chi}\left(0,\bar{\chi}_{\PR}\right),\]
where $\mathrm{R}(\MCO_F)_\chi\Fit_\MSO\left(\left(\MSO\otimes_\MBZ \Cl\left(\MCO_F\right)\right)_\chi\right)$ and 
$\Fit_\MSO\left(\left(\MSO\otimes_\MBZ\mu(F)\right)_\chi\right) \mathrm{L}_{\MFf_\chi}\left(0,\bar{\chi}_{\PR}\right)$ are viewed as $\MSO$-submodules of $\MBQ(\mu_g)$.
\end{theo}

\noindent\textsl{Proof.}
We keep the notation $\MCF_\chi := \Fit_\MSO\left(\left(\MSO\otimes_\MBZ \Cl\left(\MCO_F\right)\right)_\chi\right)$.
From Theorem \ref{theosupergras}, (\ref{bonobo}), and the property $\left.v\right)$ of index-modules, we have 
\begin{equation}
\MCF_\chi = \Fit_\MSO\left(\left(\MSO\otimes_\MBZ\left(\MCO_F^\times/\MCE_F\right)\right)_\chi\right) 
= \left[\left(\MSO\ell_F\left(\MCO_F^\times\right)\right)_\chi : \left(\MSO\ell_F\left(\MCE_F\right)\right)_\chi \right]_\MSO.
\label{spergrun}
\end{equation}
By the property $\left.iv\right)$ of index-modules, we deduce from (\ref{spergrun}) the decomposition
\begin{equation}
\mathrm{R}(\MCO_F)_\chi\MCF_\chi = \left[\MSO e_\chi : \left(\MSO\ell_F\left(\GGW\right)\right)_\chi\right]_\MSO 
\left[\left(\MSO\ell_F\left(\GGW\right)\right)_\chi : \left(\MSO\ell_F\left(\MCE_F\right)\right)_\chi\right]_\MSO.
\label{spergrdeux}
\end{equation}
By (\ref{proptrun}), and the property $\left.iii\right)$ of index-modules, we have
\begin{equation}
\left[\MSO e_\chi : \left(\MSO\ell_F\left(\GGW\right)\right)_\chi\right]_\MSO = \MSO w_\infty \mathrm{L}_{\MFf_\chi}\left(0,\bar{\chi}_{\PR}\right).
\label{OLOmega}
\end{equation}
From (\ref{spergrdeux}), (\ref{OLOmega}), and Lemma \ref{lmsuperanclfor}, we obtain
\[\mathrm{R}(\MCO_F)_\chi\Fit_\MSO\left(\left(\MSO\otimes_\MBZ \Cl\left(\MCO_F\right)\right)_\chi\right) = 
\mathrm{L}_{\MFf_\chi}\left(0,\bar{\chi}_{\PR}\right) \Fit_\MSO\left(\left(\MSO\otimes_\MBZ\mu(F)\right)_\chi\right).\]
\hfill$\square$\\

\bibliographystyle{amsplain}

\end{document}